\documentclass[12pt, A4]{article}
\usepackage{amssymb}
\usepackage{amsmath,amsfonts,amsthm,amstext}
\usepackage[english]{babel}
\usepackage[latin1]{inputenc}
\usepackage{times}
\usepackage[all]{xy}
\usepackage[symbol]{footmisc}
\DefineFNsymbols{stars}[text]{* {**} {*} {\S} {\I}}
\setfnsymbol{stars}

\newtheorem{thm}{Theorem}[section]
\newtheorem{cor}[thm]{Corollary}
\newtheorem{pro}[thm]{Proposition}
\newtheorem{deff}[thm]{Definition}
\newtheorem{lem}[thm]{Lemma}
\newtheorem{rem}[thm]{Remark}

\newcommand{\nc}{\newcommand}

\nc{\cc}{\D{C}} \nc{\hh}{\D{H}} \nc{\nn}{\D{N}} \nc{\oo}{\D{O}}
\nc{\qq}{\D{Q}}
 \nc{\rr}{\D{R}}
\nc{\zz}{\D{Z}} \nc{\livre}{\ast}

\nc{\barr}{\begin{array}} \nc{\earr}{\end{array}}
\nc{\bthm}{\begin{thm}} \nc{\ethm}{\end{thm}}
\nc{\bpro}{\begin{pro}} \nc{\epro}{\end{pro}}
\nc{\blem}{\begin{lem}} \nc{\elem}{\end{lem}}
\nc{\bins}{\begin{ins}} \nc{\eins}{\end{ins}}
\nc{\bcor}{\begin{cor}} \nc{\ecor}{\end{cor}}
\nc{\brem}{\begin{rem}} \nc{\erem}{\end{rem}}
\nc{\bdeff}{\begin{deff}} \nc{\edeff}{\end{deff}}
\nc{\bea}{\begin{eqnarray}} \nc{\eea}{\end{eqnarray}}
\nc{\D}[1]{{\mathbb#1}}
\def\R{\rm I\kern -.2em R}
\def\N{\rm I\kern -.18em N}
\def\Z{\rm Z\kern -.332em Z}
\def\de{\rm [\kern -.15em [}
\def\dd{\rm ]\kern -.15em ]}
\def\||{\hspace{0.15cm}|\hspace{0.15cm}}

\title{The Figure Eight Knot Group is Conjugacy Separable}
\author{S. C. Chagas,
 \,\,\, P. A. Zalesskii
 \footnote{\vspace*{-.5cm} Both authors were supported by
CNPq and Capes.}
  }

\begin{document}
\maketitle

\abstract{ We prove that torsion free subgroups of
$PGL_2(\mathbb{C})$ (abstractly) commensurable with the Euclidean
Bianchi groups are conjugacy separable. As a consequence we deduce
the result stated in the title. }
\section{Introduction}
\renewcommand{\thefootnote}{\arabic{footnote}}
 \setcounter{footnote}{1}

  A group $G$ is conjugacy separable if whenever $x$ and $y$ are
  non-conjugate elements of $G$, there exists some finite
  quotient of $G$ in which the images of $x$ and $y$ are
  non-conjugate. The notion of the conjugacy separability owes its importance
  to the fact, first pointed out by Mal'cev \cite{M}, that the conjugacy problem has
  a positive solution in finitely presented conjugacy separable groups.
  It is known that free products of conjugacy
  separable groups are again conjugacy separable (see \cite{R-71}). The property is
  not preserved in general by the formation of free products with
  amalgamation and $HNN$-extensions. The most general result on
  conjugacy separability of free products with cyclic amalgamation
  was proved in \cite{R-S-Z}. To investigate the property for a
  non-cyclic amalgamation is much more difficult and only few results
  exist on the subject (see \cite{R-T-V-95} and \cite{Z-T}). One of the most interesting examples of
  free amalgamated products are Bianchi groups, i.e., the groups $PSL_2(O_d)$, where
  $O_d$ denotes the ring of integers of the field $\mathbb{Q}(\sqrt{-d})$ for
  each square-free positive integer $d$. These groups have long been
  of interest, not only because of their intrinsic interest as abstract groups,
  but also because they arise naturally in number theory and geometry. In
  \cite{W-Z} conjugacy separability of Bianchi groups $PSL_2(O_d)$
  for $d=1,2,7,11$ was proved and  conjectured for the rest of them.

  In this article we prove the conjugacy separability for torsion free
  subgroups of $GL_2(\mathbb{C})$ (abstractly)
  commensurable with the Euclidean Bianchi groups, i.e., groups $PSL_2(O_d)$, $d=1,2,3,7,11$. We observe that
  it is not known whether a group  commensurable with a conjugacy separable group $G$
  (i.e., having a subgroup of finite index isomorphic with some subgroup
  of finite index of $G$) is conjugacy separable.

 As a consequence of our main result we deduce that the figure eight knot group is
 conjugacy separable. Conjugacy separability for knot groups is not known in
 general, however is known for the trefoiled knot group since it is
 a free product of cyclic groups with cyclic amalgamation.

 From a group theoretic point of view conjugacy separability of  a group
 indicates  that it has many subgroups of finite index. The other
 property that indicates this  is subgroup separability. Recall that a
  group is called subgroup separable (or LERF) if
 its finitely generated subgroups are closed in the profinite
 topology.  It follows from the definition that commensurability
 preserves  subgroup separability. Although  conjugacy separability and  subgroup separability
 do not imply each other, quite often they hold together.
   The subgroup separability of Bianchi groups was proved recently (see
 \cite[Theorem 3.6.1]{L-R}). Combining this with the results of the paper we deduce that a torsion-free group commensurable with
any finitely generated subgroup of an Euclidean Bianchi group is
conjugacy separable.

 The methods used in the paper are based on Bass-Serre theory of
 groups acting on trees and its profinite version. The basic notions
 of these theories are freely used here and can be found
 in \cite{Serre} and
 \cite{Z-M-89}, \cite{Z-89}, \cite{Z-M-90} respectively.

 \medskip
 For a group $G$ and elements $g,h\in G$ we shall frequently  use notation
 $g^h=h^{-1}gh$ in the paper.

\section{Preliminary results}


{\large\bf 2.1 Centralizers}
\begin{lem}\label{lema0} If $g \in G = GL_2(\cc), SL_2(\cc)$
 such that $g\notin Z(G)$, then $C = C_G(g)$ is abelian.
\end{lem}
\begin{proof} Let $g = \left(\begin{array}{cc}
            a & b\\
            c & d
           \end{array}\right)  \in GL_2(\cc)$, $ad-bc\ne 0$. Using the  Jordan canonical
           form we may
assume that $g=$ $\left(\begin{array}{cc}
a & b\\
0 & d
\end{array}\right),$ $ad \ne 0$. To calculate  the centralizer of $g$,
we have two cases to consider.
\vspace{0.3cm}\\
Case 1. The eigenvalues  of $g$ are distinct. Hence $g$ is
diagonal and so $g=\left(\begin{array}{cc}
            a & 0\\
            0 & d
           \end{array}\right)$.
 Then a straightforward calculation shows that the centralizer of $g$ is equal to
$$C_{GL_2(\cc)}(g)=
\left\{ \left. \left(
            \begin{array}{cc}
            x& 0\\
            0& w
           \end{array}
\right) \right| xw\ne 0 \right\}.$$
\vspace{0.3cm}\\
Case 2. If  $a = d$, i.e.,  the eigenvalues are equal, a
straightforward calculation shows that the centralizer of $g$ is
$$C_{GL_2(\cc)}=
\left\{ \left. \left(
            \begin{array}{cc}
            x& y\\
            0& x
           \end{array}
\right) \right| x\ne 0 \right\}.$$

Thus in both cases   the centralizer of a non-central element
$g\in GL_2(\cc)$ is abelian.
\vspace{0.3cm}\\
Finally observe that $C_{GL_2(\cc)}(g)\cap SL_2(\cc)=
C_{SL_2(\cc)}(g)$ and therefore the centralizer of a non-central
element of $g\in SL_2(\cc)$ is also abelian.
\end{proof}
\begin{lem}\label{lema00}
 Let $\overline{g}\in  PGL_2(\cc)$
 be an element of infinite order. Then the centralizer $C_{PGL_2(\cc)}(\overline{g})$ of $\overline{g}$ is
 abelian. Moreover, $C_{\Gamma_d}(\overline{g})$ is virtually free abelian of rank at most
 $2$, where $\Gamma_d= PSL_2(O_d)$.
 \end{lem}
 \begin{proof}
 Consider the projection  $\varphi: GL_2(\cc)\longrightarrow
PGL_2(\cc) = GL_2(\cc)/Z$, where $ Z = Z(GL_2(\cc)) = \left\{
\left. \left(
            \begin{array}{cc}
            a& 0\\
            0& a
           \end{array}
\right) \,\,\right|\,\, a\in \cc \right\}$  is the center of
$GL_2(\cc)$. We use $\overline g$ for the image of $g$ in
$GL_2(\cc)/Z$. Pick $1\neq\overline{g}\in PGL_2(\cc)$.

 We write down the inverse image of the centralizer of
the element $\overline{g}\in PGL_2(\cc)$:
$$\varphi^{-1}(C_{PGL_2(\cc)}(\overline{g}))= \{ h\in GL_2(\cc)\,\,|\,\,
hgh^{-1} = gz, \,\,\,\mbox{for some}\,\,\,z\in Z\}.$$
 Using the  Jordan canonical form we may suppose that  $g= \left(
            \begin{array}{cc}
            x& y\\
            0& w
           \end{array}
\right)$ and $ gz=\left(
            \begin{array}{cc}
            ax& ay\\
            0& aw
           \end{array}
\right)$. Since conjugate matrices have the same eigenvalues,
theses matrices  can be conjugate in $GL_2(\cc)$ only if either
$a= 1$ or  $aw =x$ and $w = ax$. In the second case we have
$a^2=1$ and so $a=-1$. But the matrix $ g=\left(
            \begin{array}{cc}
            x& y\\
            0& -x
           \end{array}
\right)$ has finite order in $PGL_2(\cc)$, because $g^2=\left(
            \begin{array}{cc}
            x^2& 0\\
            0  & x^2
           \end{array}
\right)$ is in the center of $GL_2(\cc)$. Since $g$ has infinite
order by hypothesis,  $a=1$.

Thus we have
  $$\varphi^{-1}(C_{PGL_2(\cc)}(\overline{g}))= C_{GL_2(\cc)}(g).$$
Since $Z\leq C_{GL_2(\cc)}(g)$, we have that $C_{PGL_2(\cc)}(g)$
is the quotient group of the centralizer of $g$ in $GL_2(\cc)$
modulo the center, i.e., $C_{PGL_2(\cc)}(g) \cong
C_{GL_2(\cc)}(g)/Z$. By Lemma \ref{lema0} the centralizer of $g$
in $PGL_2(\cc)$ is abelian.

The second part of the lemma follows from the fact that the
torsion free abelian subgroups in $\Gamma_d$ are free abelian  of
rank at most $2$, see Fine \cite{Fine} in page 107.
 \end{proof}

 \begin{cor}\label{diedral} A subgroup $G$ of $PGL_2(\mathbb{C})$  does not contain
 the generalized dihedral group. In particular, this is so for Bianchi groups.\end{cor}

 \begin{proof} If the generalized dihedral group  $E\cong \langle x, y | x^2 = y^2\rangle$
 is a subgroup of $PGL_2(\mathbb{C})$ then $E\leq C_{G}(x^2)$ which is an abelian group by Lemma
\ref{lema0}.\end{proof}


\subsection{ The Profinite topology}

The profinite topology on a group $G$ is the topology where the
collection of all finite index normal subgroups of $G$ serves as a
fundamental system of neighborhoods of the identity element $1\in
G$, turning $G$ into a topological group.  The completion
$\widehat{G}$ of $G$ with respect to this topology is called the
profinite completion of $G$ and can be expressed as an inverse
limit
$$\widehat{G} =\lim\limits_{\displaystyle\longleftarrow\atop N\in
{\large{\aleph}}}G/N, $$
 where $\mathcal{N} = \{ \, N \,|\, N \unlhd_f G\}$.
 Thus $\widehat{G}$ becomes a profinite group, i.e. a compact
 totally disconnected topological group. Moreover, there exists a
 natural homomorphism $\iota: G \longrightarrow \widehat{G}$  that
 sends $g \mapsto(gN)$, this homomorphism  is a monomorphism when $G$ is
residually finite. If $S$ is a subset of a topological group
$\widehat{G}$, we denote by $\overline{S}$ its closure in
$\widehat{G}$. The profinite topology on $G$ is induced by the
topology of $\widehat{G}$. Note that for a subgroup $H$ of $G$,
the profinite topology  of $H$ can be stronger than the topology
induced by the profinite topology of $G$.

Next proposition expresses the conjugacy separability property of
$G$ in terms of the profinite topology and we shall use it freely
in the paper.
\begin{pro}
 Let $G$ be a group, then the following conditions are equivalent:
 \begin{itemize}
 \item[(i)] $G$ is conjugacy separable;
  \item[(ii)] for each $x\in G$, the conjugacy class   $x^{G}$ of
  $x$ is closed in the profinite topology. In particular $G$ is residually finite;
  \item[(iii)] $G$ is residually finite and for each pair of elements $x,
y\in G$ such that $y = x^{\gamma}$,  for some $\gamma\in
\widehat{G}$, there exists $g\in G$ such that $y = x^{g}$.
 \end{itemize}
\end{pro}

 \begin{lem}\label{inter}
Let $G$ be a residually finite group, $U$ a finite index  subgroup
of $G$ and $H$ a subgroup of $G$. Then $\bar{U}\cap \bar{H} =
\overline{U\cap H}$.
\end{lem}
\begin{proof}
 We just need to proof that $\bar{U}\cap \bar{H} \leq
\overline{U\cap H}$, since clearly  $ \overline{U\cap H}\leq
\bar{U}\cap \bar{H}$.

Consider the following commutative diagram

$$
    \xymatrix{
G\ar[rr]\ar[rdd]&            &\widehat{G}\ar[ddl]\\
H\ar[rr]\ar[rd]_{\varphi_{H}}\ar[u]  &            & \bar{H}\ar[dl]^{\varphi_{\bar{H}}}\ar[u]\\
                                  & G/U= \widehat{G}/\bar{U} &}
$$
Note that ${\rm{Im}}\varphi_{\bar{H}}= {\rm{Im}}\varphi_{H}$,
${\rm{Ker}}(\varphi_{\bar{H}}) = \bar{H}\cap \bar{U}$ and
${\rm{Ker}}(\varphi_{H}) = H\cap U$. Therefore,  $[\bar{H}:
(\bar{H}\cap \bar{U}) ]= |{\rm{Im}}\varphi_{\bar{H}}|=
|{\rm{Im}}\varphi_{H}| = [H : (H\cap U)] \geq [\bar{H}:
\overline{H\cap U}]$. It follows that  $\bar{U}\cap \bar{H} =
\overline{U\cap H}$ as needed.
\end{proof}

Let $G =\pi_1({\cal{G}}, {\Gamma})$  be a fundamental group of the
  graph of conjugacy separable groups $({\cal{G}}, {\Gamma})$.
  Then the profinite topology on $G$ might induce a
weaker topology on the vertex and edge groups than their own
profinite topologies. This prevents to use  effectively the
conjugacy separability of the vertex and edge groups for proving
the conjugacy separability of $G$. Thus it is natural to assume
that
 the profinite topology on $G$ is reasonably ``strong''.
 We formalize it in  the following definition.

\begin{deff}
  Let $G$ be a group acting on a tree $S$ such that $S/G$ is finite.  We say that
  the profinite topology on $G$ is efficient, if $G$ is residually finite, the edge stabilizers $G_e$ and the vertex  stabilizers
   $G_v$ are closed in the profinite topology on $G$  and the profinite topology of $G$ induces the
  full profinite topology  on the vertex and edge stabilizers.
    \end{deff}

If $S/G$ is finite, then we can define a profinite graph
$$\widehat S=\lim\limits_{\displaystyle\longleftarrow\atop
U\triangleleft_f G} S/U,$$ where $U$ ranges over the normal
subgroups of finite index of $G$. The profinite completion
$\widehat G$ acts continuously on $\widehat S$ and $\widehat
S/\widehat G=S/G$. Moreover,   $\widehat S$ is a profinite tree
(see Proposition 3.8 in \cite{Z-M-89}).

    \begin{rem}\label{S e denso em S}
  Suppose the profinite topology on $G$ is efficient.
  Then $S$ embeds
naturally in $\widehat{S}$. This follows from the fact that
$G/G_m$  embeds in $\widehat{G}/\widehat{G}_m$ because $G_m$ are
closed in $G$ for all $m\in S$. Moreover, $S$ is dense in
$\widehat{S}$.
 \end{rem}

  We define the standard trees $S(G)$ on which $G$ acts  (respectively, $S(\widehat{G})$ on which $\widehat{G}$
  acts) only for the cases of an amalgamated free product $G=G_1\ast_HG_2$ (respectively, $\widehat G=\widehat G_1\amalg_{\widehat{H}}\widehat{G_2}$)
  and an HNN-extension $G=HNN(G_1, H,t)$ (respectively, $\widehat G=HNN(\widehat G_1,\widehat H,t)$) since we shall use them only for these cases.

  \begin{itemize}
  \item Let $G=G_1\ast_HG_2$. Then
 the vertex set is $V(S(G))= \displaystyle G/G_1\cup G/G_2$,
  the edge set is $E(S(G))= G/H$, and
  the initial and terminal vertices of an edge $gH$ are
  respectively  $gG_1$ and $gG_2$.
   \item Similarly, let $\widehat G=\widehat G_1\amalg_{\widehat{H}}\widehat{G_2}$. Then the vertex set is $V(S(\widehat{G}))=
  \displaystyle \widehat G/\widehat G_1 \cup\widehat{G}/\widehat{G}_2$,
  the edge set is $E(S(\widehat{G}))= \widehat{G}/\widehat{H}$, and
  the initial and terminal vertices of an edge $ g\widehat{H}$ are
  respectively  $g\widehat{G}_1$ and $g\widehat{G}_2$.
    \item Let $G=HNN(G_1, H,t)$. Then
 the vertex set is $V(S(G))= \displaystyle G/G_1$,
  the edge set is $E(S(G))= G/H$, and
  the initial and terminal vertices of an edge $gH$ are
  respectively  $gG_1$ and $gtG_1$.
  \item Similarly Let $\widehat G=HNN(\widehat G_1, \widehat H,t)$. Then
 the vertex set is $V(S(\widehat G))= \displaystyle \widehat G/\widehat G_1$,
  the edge set is $E(S(\widehat G))= \widehat G/\widehat H$, and
  the initial and terminal vertices of an edge $g\widehat H$ are
  respectively  $g\widehat G_1$ and $gt\widehat G_1$.
  \end{itemize}

It follows that $S(\widehat G)=\widehat{S(G)}$ and so Remark
\ref{S e denso em S}  applies for $S(G)$.

\blem\label{isolado}
  Let $G= G_1\coprod G_2\coprod \dots \coprod G_k$ be a free
  profinite product of profinite groups. Then $G_i$ are isolated in
  $G$, i.e., if $g^n\in G_i$, for some natural number $n$ and some $g\in G$ of infinite order, then $g\in G_i$.
  \elem
 \begin{proof} By an obvious induction on the number of free factors it suffices to consider the case $k=2$. Let
  $g\in G$ be an element such that $g^n
\in G_i$, we show that $g\in G_i$, for some $i=1,2$. Consider the
action of $G$ on the standard profinite tree $S(G)$ (one can use
its construction above for amalgamated free product and take an
amalgamating subgroup to be trivial). Then $g^n$ fixes a vertex in
$S(G)$. Hence $\langle g \rangle/\langle g^n\rangle$ is a finite
group acting on the subtree of fixed points $T^{g^n}$ (see Theorem
2.8 in \cite{Z-M-89}). By Theorem 2.10 in \cite{Z-M-89} $\langle g
\rangle/\langle g^n\rangle$ fixes a vertex in $T^{g^n}$ and so $g$
fixes the same vertex, say $hG_j$. It follows that $hgh^{-1}\in
G_j$. Since $G_i\cap hG_jh^{-1}=1$ for $i\neq j$ as one can easily
see mapping $G$ onto $G_1\times G_2$ we deduce that $i=j$. Since
$G_i\cap G_i^x=1$ for  $x\in G\setminus G_i$  (see Theorem 9.1.12
in \cite{RZ}) $g\in G_i$, as needed.
\end{proof}

 We distinguish two types of elements of a group $G$ acting on a tree $S$: we say that a
 non-trivial
 element $g$ is {\it hyperbolic} if it does not stabilize any  vertex of $S$. Otherwise,
 $g$ is called {\it non-hyperbolic}. The following
proposition will be used in some of our proofs.

 \begin{pro}[Proposition 24 in \cite{Serre}]
Let $G$ be an amalgamated free product or an $HNN$- extension and
$S=S(G)$ its standard tree. Suppose that an element $a\in G$ is
hyperbolic. Put $m = \displaystyle{{\rm{min}}_{v\in V(S)}{l}}\,[v,
av]$ and $T_a = \{ v\in V(S)\,\,|\,\, l[v,av] = m\}$, where
$l[v,av]$ represents the length of the geodesic $[v, av]$. Then
$T_a $ is the vertex set of a straight line (that is, a doubly
infinite chain of $S(G)$), on which $a$ acts as a translation of
amplitude $m$. Moreover, all $\langle a\rangle$-invariant subtrees
of $S(G)$ contain $T_a$. Finally if $v\in T_a$, then $T_a =
\langle a\rangle[v,av[$.
 \end{pro}

The followings proposition gives information on the closure of a
straight line in $S(\widehat{G})$, and shows that a hyperbolic
element of $G$ is hyperbolic in the profinite sense.

\begin{pro}\label{sobre retas}
 Let $G$ be an amalgamated free product or an $HNN$- extension and
$S=S(G)$ its standard tree. Suppose also that the profinite
topology on $G$ is efficient. Let $a\in G$ be a  hyperbolic
element and $T_a$ the corresponding straight line. Then:
 \begin{itemize}
 \item[(i)] $\overline{\langle a\rangle}$ acts freely on the tree
 $S(\widehat{G})$;
 \item[(ii)] $\overline{T_a} = \overline{\langle a\rangle}[v, av[$, where $v$ is a vertex of $T_{a}$;
 \item[(iii)] $\overline{\langle a\rangle}$ acts freely on the profinite tree
  $\overline{T_a}$, and $\overline{\langle a\rangle} \cong
 \widehat{\zz}$;
 \item[(iv)] $T_a$ is a connected component of \,$\overline{T_a}$
 considered as an abstract graph, in others words, the only vertices of \,$\overline{T_a}$
 that are at a finite distance from a vertex of\, $T_a$ are those of\, $T_a$.
\end{itemize}
 \end{pro}
\begin{proof}
 Item (i) is proved in  Proposition 2.9
 \cite{R-Z}. Items (ii) and (iii) are proved in  Lemma 4.1
  \cite{R-S-Z}.

 To prove (iv) we use the same argument as in  Lema~4.3
 \cite{R-S-Z}. Let $w$ be a vertex of
 $\overline{T_a}$ that is at finite distance from $v$. We show
 that $w\in T_a$. Suppose that $w \notin T_a$, then there exists an edge
  $e'$ of the geodesic $[v, w]$ that is not  in  $T_a$ and we can assume
  that the initial vertex of $e'$ is $v$. Since $\overline{T_a}=
 \overline{\langle a\rangle}[v, av[$, there exists $\alpha\in \overline{\langle
 a\rangle}$ such that $\alpha e = e'$, where $e$ is an edge of $[v,
 av[$. Let $w_0$ be the origin  of the edge $e$, then $\alpha w_0 = v$. Thus
 $\alpha w_0 = v\in T_a$, and since $T_a= \langle a\rangle[v, av[$
 there exists $\beta\in \langle
 a\rangle$ such that $\beta\alpha w_0\in[v,
 av[$, so $\beta\alpha=1$, i.e. $\alpha\in \langle
 a\rangle$. Then $\alpha e = e'\in T_a$ and so $w\in T_a$, a contradiction.
 \end{proof}

 We shall need the following generalization of the M. Hall Theorem
 to free products. Remind that a group $G$ is called LERF if  each finitely generated subgroup is closed in the
 profinite topology of $G$.


\begin{pro}[\cite{G-M-S}]\label{decomposicao compativel} Let $G_1, \dots, G_m$ be LERF groups
and let $H$ be a finitely generated subgroup of the free product
$G= G_1\ast \dots \ast G_m$. Then there exists a subgroup $U$ of
$G$ of finite index and (Kurosh-type) decompositions
$$U= U_1\ast \dots \ast U_t \,\, \mbox{and}\,\, H = H_1\ast \dots \ast H_t$$
such that
 \begin{itemize}
 \item[(a)] $H_i\leq U_i$, $i=1,\dots, t$;
 \item[(b)] For each $i=1, \dots, t-1$, $U_i$ is a subgroup
 of finite index of a conjugate of some $G_j$, $j=1, \dots, m$, i.e., $U_i= \tau G_j
 \tau^{-1}$ for some $\tau \in G$;
 \item[(c)] $U_t$ is a free group of finite rank and $H_t$ is a
 free factor of\,\, $U_t$.
 \end{itemize}
\end{pro}

Sometimes we need to separate  the conjugacy class  of a subset or
a subgroup of $G$ in a finite quotient rather than the conjugacy
class of an element. The following definition formalizes this.

 \bdeff
 A subgroup $H$ of a group $G$ is called conjugacy distinguished,
 if\,\, $\displaystyle\cup_{g\in G}H^{g}$  is closed in the profinite
 topology of $G$. Or equivalently, for every element $g\in G$ such that $g^{\gamma}\in \overline{H}$, where
 $\gamma\in \widehat{G}$, there exists $\delta\in G$ such that $g^{\delta}\in H$.
 \edeff

The following proposition collects important results about free
product of virtually abelian and virtually free groups that will
be used in the paper.
 \bpro\label{distinguido sob conjugacao}
 Let $G$ be a finitely generated free product of virtually abelian groups
 and a virtually free groups, then:
 \begin{enumerate}
\item $G$ is conjugacy separable;
 \item $G$ is LERF;
  \item The profinite topology on $G$ induces the (full) profinite topology on
  every finitely generated subgroup of $G$;
 \item For each pair of finitely generated subgroups $H_1, H_2$  of $G$,
 the product $H_1H_2$ is a closed subset of $G$, i.e., $\overline{H_1H_2}\cap G = H_1H_2$,
 where $\overline{H_1H_2}$ is the closure of $H_1H_2$ in $\widehat{G}$;
 \item Every finitely generated subgroup  $H$ of $G$ is conjugacy
 distinguished;

 \end{enumerate}
  \epro
 \begin{proof}
  Item $1$ follows from the fact that a free product of conjugacy separable groups
  is conjugacy separable and that virtually abelian and virtually free groups are
  conjugacy separable (see \cite{Dyer-79} and \cite{R-71}).
  Item $2$ follows from the fact that a free product of LERF groups is LERF
  (see \cite{Ro} and \cite{B}). Item $3$ follows from $2$, because every subgroup $H$ of
  finite index of a finitely generated subgroup $K$ is finitely generated and so is closed in the profinite topology of
  $G$. Indeed, this means that $H$ is the intersection of
  subgroups of finite index of $G$ and so for some of them, say
  $U$ one has $U\cap K=H$.
  Item $4$ was proved in \cite{R-Z-04} and \cite{C-2001}.

  $5$.  Let $g\in G$ and suppose that
$g^{\gamma}\in\bar{H}$, for some $\gamma \in \widehat{G}$.
  Let $\langle g\rangle$
 be the subgroup generated by $g$. Then by Proposition \ref{decomposicao
compativel} there is a subgroup $U$ of finite index such that
$\langle g\rangle\leq U \leq G$, and $\langle g\rangle$ is a
subgroup of a free factor of the Kurosh  decomposition for $U$.
Hence $g$ is a non-hyperbolic element of $U$. Since $U\leq_f G$,
$\widehat{G}= \widehat{U}G$, so $\gamma= u \sigma$, where $u\in
\widehat{U}$ and $\sigma\in G$. Therefore substituting $H$ by
$H^{\sigma^{-1}}$ we may assume that $\gamma\in\widehat{U}$.

Now observe that $H\cap U$ is a finitely generated subgroup of
$U$, because $H\cap U$ has finite index in $H$. Hence by
Proposition \ref{decomposicao compativel}, there exists a finite
index subgroup $U_0$ of $U$ such that $H\cap U \leq U_0 \leq U$
and the Kurosh decomposition of $H\cap U$ is compatible with the
Kurosh decomposition of $U_0$. Namely the decompositions of $U_0$
and $H\cap U$, are respectively as follows
 \begin{eqnarray*}
 U_0 &=& U_0^1\ast \ldots \ast U_0^m\\
 H\cap U&=& H_1\ast \ldots \ast H_m,
 \end{eqnarray*}
where $H_i\leq U_0^i.$
 Since $U_0\leq_f U$, $\widehat{U}=
\widehat{U_0}U$, so $\gamma = g'\delta$, where $g'\in
\widehat{U_0}$ and $\delta\in U$. Thus substituting  $U_0$, $H$
and $H\cap U$ (together with their Kurosh decompositions) by
conjugates $U_0^{\delta{-1}}$, $H^{\delta^{-1}}$ and $(H\cap
U)^{\delta^{-1}}$ in $G$, we can suppose that $\gamma\in
\widehat{U_0}$.
 Since $g^{\gamma}\in\bar{H} $ and
$g^{\gamma}\in\bar{U}$, one has $g^{\gamma}\in\bar{H}\cap\bar{U}$.

Since $U_0 \leq_f U$, then $g^n\in U_0$ for some natural number
$n$. Since $g$ is an element of a free factor of $U$, $g^n$ is an
element of a conjugate of a free factor of $U_0$, so without loss
of generality we suppose that $g^n\in U_0^1$. Now
$$\bar{U_0}= \widehat{U}_0 = \widehat{U}_0^1\amalg \ldots \amalg
\widehat{U}_0^m.$$

Then by the projective limit argument one deduces that
$g^{n\gamma}$ is conjugate to an element of a free factor of
$$\overline{H\cap U}= \bar{H_1}\amalg \ldots \amalg \bar{H_m}.$$
Indeed, $\bar U_0=\lim\limits_{\displaystyle \longleftarrow} \bar
U_{0V}$, where
$$\bar{U_{0V}}= \widehat{U}_0^1/(\widehat{U}_0^1\cap V)\amalg \ldots
\amalg \widehat{U}_0^m/(\widehat{U}_0^m\cap V)$$ and $V$ ranges
over all open subgroups of $\bar U_0$. Then $$\overline{(H \cap
U_{0})_V}= \widehat{H}_0^1/(\widehat{H}_0^1\cap V)\amalg \ldots
\amalg \widehat{H}_0^m/(\widehat{H}_0^m\cap V).$$ Since the image
$g^{n\gamma}_V$ of $g^{n\gamma}$ in $\overline{(H \cap U_{0})_V}$
has finite order, by Theorem 1 in \cite{H-R} it is conjugate in
$\overline{(H \cap U_{0})_V}$ to an element of a free factor for
every $V$. Therefore $g^{n\gamma}$ is conjugate in
$\overline{H\cap U}$ to an element of $\bar H_i$. By Theorem B' in
\cite{H-R} $\bar U_{0}^1\cap H_j^u=1$ for all $u\in \bar U_0$,
$j>1$. Hence $g^{n\gamma}$ is conjugate to $\bar H_1$ in
$\overline{U_0\cap H}$.

    Thus we may assume that $g^{n\gamma}\in
\bar{H_1}$. By Lemma \ref{inter} $\overline H\cap \overline
U_0=\overline{H\cap U_0}$,  therefore, $g^{\gamma}\in \bar{H_1}$,
since $\bar{H_1}$ is isolated in $\overline{H\cap U}=
\bar{H_1}\amalg \ldots \amalg \bar{H_m}$ (see Lemma
\ref{isolado}). Since  $g^n\in U_0^1$ and $g^{n\gamma}\in U_0^1$,
by Theorem 9.1.12 in \cite{R-Z} the conjugating element $\gamma\in
\bar{U_0^1}$. It follows that $g\in U_0^1=G\cap \overline{U_0^1}$.
 But $U_0^1$ is either virtually free or virtually abelian, so $H_1$ is conjugacy distinguished in
$U_0^1$. This finishes the proof.
  \end{proof}

If $\Gamma$ is a profinite graph then its connected component is a
maximal connected profinite subgraph of $\Gamma$. Note that it
coincide with the closure of a connected component of $\Gamma$
considered as an abstract graph, since the closure of a connected
abstract subgraph is a connected profinite subgraph. If a
profinite group $G$ acts on $\Gamma$ and $C$ its connected
component we denote by $Stab_G(C)$ the maximal subgroup of $G$
leaving $C$ invariant.

\begin{lem}\label{component}
 Let $G$ be a profinite group acting on a profinite graph $S$ and let $m_1, m_2$ be elements of a connected
   component
  $C$ of $S$. Then $$C/Stab_G(C)\subseteq S/G.$$
\end{lem}

\begin{proof}
 We need to show that if $g\in G$ with $gm_1 =m_2$, then $g$ leaves $C$ invariant, i.e.
  $g\in Stab_{G}(C)$. Since $m_2\in C\cap gC$, $C\cup gC$ is connected and so $C=gC$.
\end{proof}

\begin{lem}\label{componentstabilizer} Let $G$ be a  group  acting  on a profinite tree
  $S$ such that $S/G$ is finite. Let  $\Delta$ be a profinite subgraph of
  $S$ and $e$ be an edge in $S$. Suppose  a connected component $C$ of $\Gamma:=\Delta \setminus Ge$ contains at least
  one edge. Then there exist edges $e'\in C$
and $e_0\in Ge\cap \Delta$ that have a common vertex $v$.
   \end{lem}

\begin{proof}
  Consider   the
  profinite graph $\Delta_{\Gamma}$ obtained from $\Delta$ by collapsing all connected
  components of $\Gamma$. By Proposition 3.2 in  \cite{Z-92}
  $\Delta_{\Gamma}$ is a profinite tree.  Let $v_{C}$ be the vertex that is the image
  of the connected component $C$. Now since
  $E(\Delta_{\Gamma})$ is compact (because $\Delta\!\setminus{\widehat{G}e}$
   is open and closed), by Proposition
   2.15 in \cite{Z-M-89} there exists an edge
  $\overline{e_0}\in \Delta_{\Gamma}$ that has $v_{C}$ as a
  vertex, say $v_{C}= d_0(\overline{e}_0)$.
  Let $e_0= \pi^{-1}_{\Gamma}(\overline{e_0})$, where
  $\pi_{\Gamma}: \Delta \longrightarrow \Delta_{\Gamma}$ is a
  quotient map. Then $d_0(e_0)\in C$ and since $E(C)= E(S(G))\cap C$ is
  compact, by Proposition 2.15 \cite{Z-M-89} there exist $e'\in E(C)$ having $d_0(e_0)$ as a
  vertex. Note also that since $\pi_{\Gamma}(E(\Delta))=\pi_{\Gamma}(\widehat Ge)$ one has $e_0\in \widehat Ge$.
   \end{proof}

 \begin{lem}\label{estabilizadores d c}
    Let $G$ be a  group that acts on a tree
  $S$ such that $S/G$ is finite and the edge stabilizers are finitely generated. Suppose that the
  profinite topology on $G$ is efficient and there exists an epimorphism
 $\tau:G\longrightarrow K$ such that
 the restriction of $\tau$ to the vertex stabilizer $G_v$ is an isomorphism for each
 $v\in V(S)$. Suppose every finitely generated subgroup of $K$ is conjugacy distinguished.
 Then the edge stabilizers of $G$ are conjugacy distinguished in $G$.
  \end{lem}
  \begin{proof}  Let $e\in S$ be  an edge. We have to show that
  $G_e$ is conjugacy distinguished in $G$. Let
  $g\in G$ and $\gamma\in \widehat{G}$ be such that
  $g^{\gamma}\in \widehat{G_e}$, we need to show that there exists
  $\gamma_1\in G$ such that $g^{\gamma_1}\in G_e$. Observe that by  item~(i)
  of Proposition \ref{sobre retas} $g$ is  a non-hyperbolic
  element of $G$, therefore $g$ stabilizes a vertex $v$ in $S$. Fix a connected transversal $T$ of $S/G$ in $S$.
  Without loss of generally
  we can suppose that $v$ and $e$ belong to $T$. We use induction
  on the number of edges of the geodesic $[v,e]$ in $T$.

  Suppose  $[v,e]$ has one edge only.
  Since $\tau$ restricted to $G_v$  is an isomorphism to $K$,
   $\tau(G_e)$ is conjugacy distinguished in $K$
   by Proposition \ref{distinguido sob
  conjugacao} ($\tau(G_e)$ is finitely generated). Hence the exists
   $k\in K$ such that $\tau(g)^{k}
  \in \tau(G_e)$ and since $\tau|_{G_v}$ is an isomorphism to K,
  there exists $\delta\in G_v$ such that $\tau(\delta)= k$. Then
  $\tau(g)^{\tau(\delta)}\in\tau(G_e)$, and consequently $g^{\delta}\in
  G_e$.

     \medskip

  Suppose now that $[v,e]$ has more than one edge.
    Since $g$ fixes $v$ and $\gamma^{-1} e$
  by Theorem 2.8 in \cite{Z-M-89} $g$ fixes the geodesic $[v, \gamma^{-1} e]$, i.e.,
   the minimal profinite subtree which contains $v$ and $\gamma^{-1} e$.
  Denote by $C$ the connected
  component of
    $[v, \gamma^{-1} e]\!\setminus\widehat{G}e$ containing $v$. Note that $C$ contains at least one edge since
    otherwise $v$ would be a vertex of a translation of $e$, contradicting the assumption that $[v,e]$ has more
    than one edge. Then by Lemma \ref{componentstabilizer}
    there are edges  $e'\in C$ and
      $e_0\in \widehat{G}e\cap [v, \gamma^{-1} e]$ that have a vertex in common, say $d_0(e_0)$.

  Let $\nu : S(\widehat{G})\rightarrow S(\widehat{G})/\widehat{G}$
  be the natural epimorphism. Observe that $\nu(e)= \widehat{G}e$ and
  $\nu(v)\in (S/G)\!\setminus\nu(e)$.  Let $\widehat\Pi$ be the maximal subgroup of
  $\widehat G$ that leaves invariant the connected component $\widehat C$ of
  $S(\widehat G)\!\setminus\widehat Ge$ containing $v$, i.e.
  $\widehat \Pi=Stab_{\widehat G}(\widehat C)$.
  By Lemma \ref{component}  $\nu|_{\widehat C}$ coincides with the factorization of $\widehat C$
   modulo $\widehat\Pi$. Therefore, replacing $\gamma$ with its
   multiple $x\gamma$ for some $x\in \widehat\Pi$ we may assume
   that $e'\in T$. But $[v,e']$ has less edges then $[v,e]$ so by
   the induction hypothesis $g$ is conjugate to an element of
   $G_{e'}$ in $G$. Hence we can assume that $g\in G_{e'}$. It
   follows that $g\in d_0(e_0)$ and replacing $v$ by $d_0(e_0)$ we
   may assume that $v$ is the initial vertex of $e_0$. But
   changing $\gamma$ again we may assume that $e_0$ is in $T$ (we
   loosing $e\in T$ of course) , so by the induction hypothesis
   again we deduce that $g$ is conjugate in $G$ to an element of
   $G_{e_0}$. However, we observe now that $T$ contains unique
   edge of the orbit $\widehat Ge$ so our new $e_0\in T$ is the
   same as our old $e$. Thus $g$ is conjugate to an element of
   $G_e$ in $G$ as required.
\end{proof}


\begin{thm}\label{G s.c.}
  Let $G$ be a torsion free group that does not contain subgroups
  isomorphic to a generalized dihedral.
  Let $H$ be a finite index subgroup of $G$ such that:
     \begin{enumerate}
  \item $H$ is conjugacy separable;
    \item for each $h\in H$, $C_{H}(h)$ is a free abelian group
    of rank at most $2$.
  \end{enumerate}
  Then $G$ is conjugacy separable.
  \end{thm}
  \begin{proof}
     Let $g_1, g_2\in G$ such that $g_2 = g_1^{\gamma}$, for some
     $\gamma\in \widehat{G}$. Since $H$ has finite index in $G$,
      $g_1^m, g_2^m \in H$, where $m=|G:H|$.

 Observe that $\widehat{G} = G\widehat{H}$, so that we can write
  $\gamma = d\gamma_0$, where $\gamma_0\in \widehat{H}$ and
$d\in G$. Therefore  $g_2^m = (g_1^{\gamma})^m =(g_1^{m})^{\gamma}
= (g_1^m)^{d\gamma_0}$. Now substituting $g_1$ by $g_1^d$, we can
suppose that $\gamma\in \widehat{H}$. Thus $g_1^m$ and $g_2^m$ are
conjugated in $\widehat{H}$, and since $H$ is conjugacy separable
there exists $h\in H$ such that $ g_1^m= (g_2^{h})^m$. Hence $
g_1^m = (g_2^m)^{h}$, so $g_1^m$ and $g_2^m$ are conjugate in $H$.
Thus we can suppose that $g_1^m = g_2^m$. Let $N = \langle g_2^m
\rangle$ and $K = \langle g_1, g_2\rangle$. Then $K$ centralizes
$N$. Since $C_{G}(N)$ is a virtually free abelian group of rank at
most $2$, $K$ is a virtually free abelian group of rank at most
$2$. If $K$ is virtually cyclic, then since $K$ is torsion free,
by Theorem 3.5 in \cite{Dicks} $K$ is cyclic and therefore $g_1 =
g_2$.

If $K$ is abelian,  $g_1= g_2$ and we one done. Suppose that $K$
is a non-abelian group having two generated free abelian subgroup
of finite index. Choose $A$ a normal torsion free abelian subgroup
of finite index, where without lost generality we can assume that
$g_1^m$ is a generator of $A$. Let $\varphi: K \longrightarrow
Aut(A)=GL_2(\zz)$ be the homomorphism induced by the action of $K$
on $A$. Since $N$ is central in $K$, $\varphi(K)$ is a  finite
subgroup of the group of the upper triangular matrices with first
column  $\left(%
\begin{array}{c}
  1 \\
  0 \\
\end{array}%
\right)$. Hence the other element of the diagonal of this matrix
must be  $1$ ou $-1$. Since the image of $K$ in $GL_2(\zz)$ is
finite and a matrix of the form
$ \left(%
\begin{array}{cc}
  1 & x\\
  0 &  1\\
\end{array}%
\right) $ has infinite order, any element of $\varphi(K)$
 has the form
 $
\left(%
\begin{array}{cc}
  1 & x\\
  0 & -1\\
\end{array}%
\right)$. But
 $
\left(%
\begin{array}{cc}
  1 & x\\
  0 & -1\\
\end{array}%
\right)
\left(%
\begin{array}{cc}
  1 & y\\
  0 & -1\\
\end{array}%
\right)
= \left(%
\begin{array}{cc}
  1 & y-x\\
  0 & 1\\
\end{array}%
\right),
 $
 so we conclude that $\varphi(K)$ is a  group of order $2$. Hence
 the centralizer $P = C_{K}(A)$ of
$A$ in $K$ has index $2$. Since $ A \leq Z(P)$,  $[P, P]$ is
finite and since $P$ is torsion free,  $P$ is abelian. Therefore
we can assume that $P = A$. In this case  conjugating the image of
$K$ in $GL_2(\zz)$ if necessary (i.e., choosing another base for
$A$) we can suppose that $\varphi(K)$ is generated either by the
diagonal matrix
$\left(%
\begin{array}{cc}
  1 & 0 \\
  0 & -1 \\
\end{array}%
\right)$ or by the matrix $\left(%
\begin{array}{cc}
  1 & 1 \\
  0 & -1 \\
\end{array}%
\right)$. Indeed, conjugating the matrix $\left(%
\begin{array}{cc}
  1 & x\\
  0 & -1\\
\end{array}%
\right)$ by the matrix $\left(%
\begin{array}{cc}
  1 & y\\
  0 & -1\\
\end{array}%
\right)$ we get $\left(%
\begin{array}{cc}
  1 & 2y-x\\
  0 & -1\\
\end{array}%
\right)$,  then we get the first matrix in the case of $x$  being
even, and the second matrix in the case of $y$  being odd; in the
latter
case  the matrix is conjugate of the matrix $\left(%
\begin{array}{cc}
  0 & 1 \\
  1 & 0 \\
\end{array}%
\right)$.

In the first case $K$ is a generalized dihedral group. Indeed,
 let $z$ be an element of the inverse image of
$\left(%
\begin{array}{cc}
  1 & 0 \\
  0 & -1 \\
\end{array}%
\right)$. Note that $Z(K)$ is cyclic, because otherwise $Z(K)$ is
of finite index and since $K$ is torsion free, $K$ would be
abelian. Therefore since $z$ centralizes $z^2$, $z^2$ is a power
of $g_1^m$. But this power can be only $1$ or $-1$, otherwise
factoring out $z^2$ we get that $K/\langle z^2\rangle$ has the
followings generators $g_1^{m}\langle z^2\rangle$, $\alpha\langle
z^2\rangle$ and $z\langle z^2\rangle$, where $\alpha$ is the
generator of the group $A$ inverted by the matrix above; hence we
have a $3$-generated group, unless $z^{2}$ is $g^m_1$ or
$(g_1^{m})^{-1}$. Thus $K$ has only one relation $\alpha^z =
\alpha^{-1}$. Hence $K$ is a generalized dihedral group which is
discarded by the hypothesis of the theorem.

In the second case, $K$ has torsion. Indeed, let $z$ be the
inverse image of the
$\left(%
\begin{array}{cc}
  0 & 1 \\
  1 & 0 \\
\end{array}%
\right)$, then there exist generators $x$ and $y$ of $A$ that are
exchanged by the action of $z$, i.e., $x^z = y$ and $y^{z}= x$.
Since $z$ centralizes $z^2$,  $z^2= (xy)^n$, for some $n\in \zz$.
But  $(zx^{-n})^2 = zx^{-n}zx^{-n} =z^2(yx)^{-n} =1$, i.e.,
$zx^{-n}$ has order $2$. Therefore $K$ has torsion and since $G$
is torsion free this case is also discarded.

Thus $K$ is isomorphic to $\zz\times \zz$, and this implies that
$g_1 = g_2$.
  \end{proof}
\section{Principal results}


Recall that the Bianchi Groups are $PSL_2(O_d)$, where $O_d$ is
the ring of integers in the  imaginary quadratic number field
$\mathbb{Q}(\sqrt{-d})$, and
$$ O_d = \mathbb{Z} + w \mathbb{Z}, \,\, \mbox{where}\,\, \begin{cases}
 w= \sqrt{-d}, \,\mbox{se}\, -d\equiv 1 \mbox{mod}\,(4)\\
 w= \frac{1+\sqrt{-d}}{2}, \,\mbox{se}\, -d\not\equiv 1
 \mbox{mod}\,(4).
 \end{cases}
$$
and $d$ is a positive square-free integer. Let $m\in \mathbb{N}$
and $O_{d,m}= \mathbb{Z} + mw\mathbb{Z}$, where $w$ is given
above. It is easy to see that $|O_d: O_{d,m}|= m$, and we can
consider the groups $SL_2(O_{d,m})$. Observe that $SL_2(O_{d,m})$
is of finite index in $SL_2(O_d)$.

\medskip

We shall use the following decompositions of Bianchi groups
$\Gamma_d = PSL_2(O_d)$ for $d=1,2,7,11$ (see \cite{Fine}).
 \begin{eqnarray}
 PSL_2(O_1) = G_1\livre_{M}G_2,\notag
 \end{eqnarray}
where $G_1= S_3\livre_{\zz/3\zz}A_4$, $G_2=
S_3\livre_{\zz/2\zz}D_2 $ and $M$ is the modular group
$PSL_2(\zz)$.

\begin{eqnarray}
PSL_2(O_2) &=& HNN\big(K_2, M, t\big)\notag\\
PSL_2(O_7) &=& HNN\big(K_7, M, t\big)\notag\\
PSL_2(O_{11}) &=& HNN\big(K_{11}, M, t\big),
 \notag
 \end{eqnarray}
where $K_2= (A_4\livre_{\zz/2\zz}D_2), K_7=
(S_3\livre_{\zz/2\zz}S_3), K_{11}= (A_4\livre_{\zz/3\zz}A_4)$ and
$M$ is the modular group $PSL_2(\zz)$.

The Bianchi group $\Gamma_3$ does not decompose in amalgamated
free product or $HNN$-extension. Therefore we are going to use the
decomposition found by Scarth of a certain subgroup of finite
index in $\Gamma_{3}$. The decomposition is described in the next

\begin{pro}[Scarth, M. R.]\label{aluno do mason} There exists a subgroup
$\Gamma_{3,2}:=PSL_2(O_{3,2})$ of $\Gamma_3$ of finite index which
 is an $HNN$-extension of $K_{3,2}$ with the modular group associated,
where $K_{3,2}= S_3\ast_{\zz/3\zz}D(3,3,3)$ and $D(3,3,3)= \langle
x, y | x^3=y^3=(xy)^3=1\rangle$.
\end{pro}
\begin{proof}
The group $PSL_2(O_{3,2})$ has the following presentation (see
\cite{Scarth})
$$\langle a, t, w | a^2 = (at)^3 = (w^{-1}awa)^3= [t,w]=1\rangle.$$

Let $v= w^{-1}a w$, then using Tietze transformations ($(av)^3=1$
iff $(va)^3=1$ and $(at)^3=(wvw^{-1}t)^3=(wvtw^{-1})^3=1$ iff
$(tv)^3=1$) we get the following presentation
$$\langle a, t, v , w | a^2= (at)^3= (av)^3= v^2 = (tv)^3, t= w^{-1}tw, v= w^{-1}aw\rangle.$$
Let $K_{3,2}= \langle a, t, v | a^2 = (at)^3 = (av)^3 = v^2
=(tv)^3\rangle$. We show that $PSL_2(O_{3,2})$ is an
$HNN$-extension of $K_{3,2}$ with $\langle a, t\rangle \cong M
\cong \langle v, t\rangle$ associated.

Observe first that $\langle a, t\rangle\cong M$, since have the same
presentation. Then $PSL_2(O_{3,2})=HNN(K_{32},M,w)$.

Let $s= at, m= av$. We get
$$K_{3,2}= \langle a, s, m | a^2= s^3= m^3= (am)^2 = (sm^2)^3=1\rangle,$$
where we used again that $sm^2=atva$ is of order $3$ iff $tv$ is
of order $3$. Observe that $\langle a, m | a^2= m^3 =
(am)^2=1\rangle \cong S_3$ and a group with presentation $\langle
s, \overline{m} | s^3= \overline{m}^3
=(s\overline{m}^{-1})^3\rangle\cong D(3,3,3).$ Then we have
$K_{3,2} = S_3\ast_{m=\overline{m}} D(3,3,3)$.
\end{proof}


\begin{lem}\label{rest. de K32 e monomorfismo}
There exists a homomorphism $\tau : \Gamma_{3,2}\rightarrow P$,
where $P$ is virtually a free product of torsion free abelian groups
and a free group such that $\tau_{|K_{32}}$ is injective.
\end{lem}
\begin{proof} We use the notation of the preceding proposition.
Let $\psi_1: S_3\rightarrow S_3$ be the map that
 sends $a\mapsto am$, $m\mapsto m^{-1}$. Observe that $\{\psi_1(a), \psi_1(m)\}$ satisfies the relations of
 the group $S_3$, so $\psi_1$ is an epimorphism  by Von Dick's theorem. Then $\psi_1$ is an
 automorphism of order 2
 because $am$ and $m^{-1}$ generate $S_3$. Let $\psi_2: D(3,3,3)\rightarrow D(3,3,3)$ be the map that
 sends $s\mapsto s^{-1}$, $\bar {m}\mapsto \bar{m}^{-1}$, again we have that $\psi_2$ is an
 automorphism of order 2. Note that these two automorphisms agree on the amalgamated subgroup $\langle m\rangle=\langle \bar m\rangle$,
 i.e., $\psi_1(m)= m^{-1}= \bar{m}^{-1}= \psi_2(\bar{m})$,
 therefore by the universal property of an amalgamated free product there exists a
 unique homomorphism $\psi:  S_3\ast_{\mathbb{Z}/3\mathbb{Z}} D(3,3,3)\rightarrow
 S_3\ast_{\mathbb{Z}/2\mathbb{Z}} D(3,3,3)$ extending $\psi_1$ and $\psi_2$. Since $\psi_1$ and $\psi_2$ are
 automorphisms  of order 2 it follows that $\psi$ is an automorphism of order 2 as well.

 Consider the group $P = (S_3\ast_{\mathbb{Z}/3\mathbb{Z}} D(3,3,3))\rtimes
 \langle\psi\rangle$. Then the map $\varphi$ that sends $K_{3,2}\mapsto
 K_{3,2}$ identically and $w\mapsto \psi$, extends to an epimorphism $\tau: \Gamma_{3,2}\rightarrow P$
 by Von Dycks theorem because $\{\varphi(a), \varphi(m), \varphi(\overline{m}), \varphi(s), \varphi(w) \}$
 satisfies the relations of $\Gamma_{3,2}$. In addition, $\tau|_{K_{3,2}}$ is a monomorphism
 to $P$. Now it remains to note that by the Kurosh subgroup theorem (see \cite{Serre} for the version with amalgamation)
 any torsion free subgroup of finite index
 of $P$ is a free product of virtually abelian group and a  free group, since $D(3,3,3)$ is virtually
 free abelian.
\end{proof}

\begin{pro}\label{intersection} Let
$H$ be finitely generated subgroup of $K_{32}$ and $M$ be the
modular group. Then $\overline{H\cap M} = \overline{H} \cap
\overline{M}$.
\end{pro}

\begin{proof}

We must prove that $\overline{H\cap M} \geq \overline{H}\cap
\overline{M}$, since, clearly $\overline{H\cap M} \leq
\overline{H}\cap \overline{M}$.

Consider  $L =\langle K_{32},K_{32}^{w^{-1}}\rangle=
K_{32}\ast_{M}K_{32}^{w^{-1}}$ be a amalgamated free product of
$K_{32}$, where $M$ is a modular group and set $P = \langle H,
K_{32}^{w^{-1}}\rangle$. By the subgroup theorem for amalgamated
free products we have $P= H\ast_{H\cap M} K_{32}^{w^{-1}}$. Since
$\Gamma_{32}$ is LERF, all the subgroups $ H, M, P$ are closed in
the profinite topology of $\Gamma_{32}$ and this topology induces
an efficient profinite topology on $P$ and on $L$.

It follows that $\widehat{L}=
\widehat{K_{32}}\amalg_{\widehat{M}}\widehat{K_{32}}^{w^{-1}}$ and
$\overline{P}= \widehat{H}\amalg_{\widehat{H\cap
M}}\widehat{K_{32}^{w^{-1}}}$.

By Exercise 9.2.7 (2) in \cite{RZ}
$\widehat{H}\amalg_{\widehat{H\cap M}} \widehat{K_{32}^{w^{-1}}}$ is
proper. Obviously, $\widehat{H}\cap \widehat{K_{32}^{w^{-1}}}$
contains $\widehat{H}\cap \widehat{M}$.
 We show that
$\widehat{H}\cap\widehat{K_{32}^{w^{-1}}}= \widehat{H\cap M}$.
 Indeed, suppose not. Pick an element
$h\in
\widehat{H}\cap\widehat{K_{32}^{w^{-1}}}\setminus\widehat{H\cap
M}$. Then the exist an open normal subgroup $U$ of $\overline{P}$
such that $hU\not\in \widehat{H\cap M}U/U$. Put $N= U\cap H$, $N'=
U\cap K_{32}^{w^{-1}}$ and $V=\widehat{H\cap M}\cap U$. Consider
the natural homomorphism $f: P \longrightarrow
H/N\amalg_{\widehat{H\cap M}/V} K_{32}^{w^{-1}}/N'$. Then
$f(h)\not\in \widehat{H\cap M}/V$. But $f(h)\in
f(\widehat{H}\cap\widehat{K_{32}^{w^{-1}}})\subset (H/N)\cap
(K_{32}^{w^{-1}}/N')= \widehat{H\cap M}/V$. So $h\in
Ker(f)\widehat{H\cap M}\leq \widehat{H\cap M}U$ contradicting with
the fact that $hU\not\in \widehat{H\cap M}U/U$. Thus
$\widehat{H}\cap\widehat{K_{32}^{w^{-1}}}= \widehat{H\cap M}$.
\end{proof}


 \begin{pro}\label{comenscomgamma}
  Let $H$ be a  finite index torsion free subgroup of $G = \Gamma
  _d$, $d = 1,2,7,11$ or $G= \Gamma_{3,2}$. Then $H$ is conjugacy separable.
  \end{pro}
  \begin{proof}
  Let $S(G)$ and $S(\widehat{G})$ be the trees associated with decompositions of
   $G$ and $\widehat{G}$.
By \cite[Theorem 3.6.1]{L-R} $G$ is subgroup separable, therefore
the profinite topology on $G$ is efficient, so $S(G)$ is embedded
in $S(\widehat{G})$ (see Remark \ref{S e denso em S}).

  By  Lemma \ref{lema00} the centralizers of non-trivial elements of $H$
are free abelian groups of rank at most $2$. Hence by Theorem
\ref{G s.c.} it is enough to proof the proposition for an
appropriate subgroup of finite index of $H$. An appropriate
subgroup  for us is a subgroup that satisfies the conditions of
Lemma \ref{estabilizadores d c}. We show  the existence of such a
subgroup. Indeed, by Lemmas 4.2 and 4.3 in \cite{W-Z} and  by
Lemma \ref{rest. de K32 e monomorfismo} in the case $d=3$, there
exists virtually a free product $P$ of torsion free abelian groups
and a free group and an epimorphism  $\tau: G \rightarrow P $ such
that $\tau|_{G_v}$ is injective, for all $v\in V(S)$ and
$|P:\tau(G_v)|< \infty$. Therefore replacing $H$ by a subgroup of
finite index if necessary we may assume that the image of $H$ is a
free product of torsion free abelian  groups and a free group. As
$S(G)/H$ is finite, there exists a finite number of orbits $Hv$,
$v\in S(G)$, i.e. there exists a finite number of vertex
stabilizers in $H$ up to conjugation. Since the index $[P :
\tau(G_v)]$ is finite, so is $[P: \tau(H_v)]$. It follows that
$\tau(H_v)$ has a finite number of conjugates in $P$. Therefore
$K:=\bigcap_{v\in V(S)}\tau(H_v)$ has finite index in $P$  and
$L=\tau^{-1}(K)\cap H$ is the desired subgroup. Indeed,
$\tau|_{L_v}$ is an isomorphism to $K$, because $\tau(L_v)=
\tau(L\cap H_v)= \tau(L)$, where the last equality follows from
the definition of $L$ (to see that the right handside is contained
in the left handside let  $l\in L$; then $\tau(l)\in K\leq
\tau(H_v)$; hence there exist $h_v\in H_v$, $x\in Ker\,(\tau)$
such that $l=xh_v$; but $Ker\,(\tau)\leq \tau^{-1}(K)$, so $h_v\in
L$). Moreover, the vertex and edges stabilizers in $L$ are
finitely generated and since $G$ is subgroup separable
\cite[Theorem 3.6.1]{L-R}, the profinite topology on $G$ induces
the  efficient profinite topology on $L$. Note also that every
finitely generated subgroup of $K$ is conjugacy distinguished by
Proposition \ref{distinguido sob conjugacao}, because $K$ is
isomorphic to a subgroup of $P$ and so by the Kurosh subgroup
theorem satisfy the hypothesis of Proposition \ref{distinguido sob
conjugacao}. Thus replacing $H$ by $L$ if necessary, from now on
we may assume that $H$ satisfies the hypotheses of
Lemma~\ref{estabilizadores d c}.

    Let $g_1, g_2\in H$ such that $g_2 = \gamma^{-1} g_1 \gamma$
for any $\gamma \in \widehat{H}$. We show that $g_1$ and $g_2$ are
conjugate in $H$. We split the proof into two cases.
   \vspace{0.3cm}\\
{\bf Case 1}:  $g_1$ is  non-hyperbolic. By item (i) of
Proposition 2.9 in \cite{R-Z} $g_2$ is  non-hyperbolic as well.
Let $v_1$ be a vertex fixed by $g_1$. If $\gamma\in
\widehat{G}_{v_1}$, then $g_2= g_1^{\gamma}\in \widehat{G}_{v_1}$.
Since $\widehat{G}_{v_1}\cap H = G_{v_1}\cap H$, the elements
$g_1, g_2\in G_{v_1}\cap H$.
 By  Lemma \ref{inter} $\widehat{G}_{v_1}\cap \widehat{H}=
\widehat{G_{v_1}\cap H}$, so  $\gamma\in \widehat{G_{v_1}\cap H}$.
Since $H$ satisfies the hypothesis of Lemma \ref{estabilizadores d
c}, the group  $H_{v_1}=G_{v_1}\cap H$ satisfies hypothesis of
Proposition \ref{distinguido sob conjugacao} and so it is
conjugacy separable. Thus there exists $\gamma_1\in G_{v_1}\cap H$
such that $g_2=g_1^{\gamma_1}$, i.e. $g_1$ and $g_2$ are conjugate
in $H\cap G_{v_1}$.

Suppose now that $\gamma\notin \widehat{G}_{v_1}$. We shall prove
then that some conjugates of $g_1$ and $g_2$ in $H$ both stabilize
some edge $e\in S$.

By Theorem 3.12 \cite{Z-M-89} $g_1= \gamma g_2\gamma^{-1}\in
\widehat{G}_{v_1}\cap\gamma \widehat{G}_v\gamma^{-1}\leq
\widehat{G}_{\bar e}$, where $\bar e\in S(\widehat{G})$, $v$ is a
vertex of $S(G)$ stabilized by $g_2$ and $v= d_0(\bar e)$ or $v=
d_1(\bar e)$. Now $|G_v~:~H_v|< \infty$, so
$\widehat{G}_v/\widehat{H}_v= G_v/H_v$ and we get
$Star_{S(\widehat{G})}(v)/\widehat{H}_v= Star_{S(G)}(v)/H_v$. This
implies the existence of an element $\widehat{h}\in \widehat{H_v}$
such that $e_1:= \widehat{h}\bar e\in S(G)$.

Then $g_1^{\widehat{h}^{-1}}\in \widehat{H_{e_{1}}}$ and since
$\widehat{h}\in \widehat{H_v}$, and $H_{e_1}$ is conjugacy
distinguished in $H_v$ (by item 5 of Proposition \ref{distinguido
sob conjugacao}), there exists $h\in H_v$ such that $g_1^{h}\in
H_{e_1}$.

 Thus we may assume that $g_1\in H_e$, where $e$ is an
edge of $S(G)$. By Lemma \ref{estabilizadores d c} the edge
stabilizers $H_e$ are conjugacy distinguished in $H$. Therefore
conjugating $g_2$ by an element of $H$ if necessary
 we may assume  that $g_2\in H_{e}$.

Recall now that $\tau: H \longrightarrow K$ is a homomorphism
satisfying hypothesis of  Lemma \ref{estabilizadores d c}.  As
$g_1, g_2$ are conjugate in $\widehat{H}$, the elements
$\tau(g_1)$ and $\tau(g_2)$ are conjugate in $\widehat{K}$, and
since $K$ is conjugacy separable there exists $k\in K$ such that
$\tau(g_2) = \tau(g_1)^k$. Then there exists $z\in H\cap G_{v_1}$
such that $\tau(z) = k$, so $\tau(g_2)=\tau(g_1^z)$, and since
$g_2, g_1^z\in G_{v_1}\cap H$ and $\tau|_{G_{v_1}}$ is injective,
one has $g_2= g_1^z$.
  \vspace{0.3cm}\\
{\bf Case 2}: The elements $g_1, g_2$ are  hyperbolic.  Let
$T_{g_1}$ and $T_{g_2}$ be the infinite straight lines on which
$g_1$ and $g_2$ act.  Let $e$ be any edge of $T_{g_1}$. Since $g_2
= \gamma^{-1} g_1\gamma$, the group $\langle g_2\rangle$ acts
freely over $\gamma^{-1} \overline{T_{g_1}}$, and by Proposition
2.2(ii) in \cite{R-Z} $\gamma^{-1}\overline{T}_{g_1} =
\overline{T}_{g_2}$. It follows that $\gamma^{-1} e \in
\overline{T}_{g_2}= \overline{\langle g_2\rangle}T_2$, where $T_2=
[v_2, g_2v_2]$ for any $v_2\in T_{g_2}$. Let $x\in
\overline{\langle g_2\rangle}$ such that $x\gamma^{-1} e\in T_2$
and $\widehat{\mu}:S( \widehat{G})\rightarrow
S(\widehat{G})/\widehat{H}$ be the natural map. Then
$\widehat\mu(x\gamma^{-1} e)=\widehat\mu(e)$, and since
$S(\widehat{G})/\widehat{H}= S(G)/H$ (as one easily checks using
the closedness of $HG_e$ in $G$)  there exists $g\in H$ such that
$x\gamma^{-1} e = g e\in S$. Therefore $g^{-1}x\gamma^{-1} e = e$
and so $ \delta:=g^{-1}x\gamma^{-1} \in \overline{H_e}$. Note that
$g_1= \gamma g_2\gamma^{-1}= \gamma x^{-1}g_2 x\gamma^{-1} =
\delta^{-1}g^{-1}g_2g\delta$. Then substituting $g_2$ with
$g^{-1}g_2g$ and $\gamma$ with $\delta$ we can suppose that
$\gamma\in \overline{H_{e}}$ and so $T_{g_1}$ and $T_{g_2}$ have a
common edge $e$.

Let $P$ be a geodesic of finite length in $T_{g_1}$ that has $e$
as one of its edges and such that $\gamma\in\overline{I}$, where
$I =\underset{e\in E(P)}{\cap}H_e$. We show that we can assume
that $\gamma$ can be substituted by  an element that belongs to
the closure of intersection of the edge stabilizers of a geodesic
which strictly contains $P$. Let $e_1\in T_{g_1}\setminus P$ be an
edge connected to $P$, and let $v$ be the common vertex of $e_1$
and $P$. Let $P^{+}= P\cup \{ e_1\}$, and put $e_2= \gamma^{-1}
e_1$ so that $e_2\in \overline{T}_{g_2}$, because
$\gamma^{-1}\overline{T}_{g_1}= \overline{T}_{g_2}$. In fact,
$e_2\in T_{g_2}$. Indeed, $\gamma^{-1}[e, e_1] = [e, e_2]$ is
finite and contained in $\overline{T_{g_2}}$. Now observe that
$T_{g_2}$ is a connected component of $\overline{T}_{g_2}$
considered as an abstract graph, in others words the only vertices
of $\overline{T}_{g_2}$ that are at a finite distance from
$T_{g_2}$ are the vertices of $T_{g_2}$, see Proposition
\ref{sobre retas} (iv). Therefore $e_2\in T_2$.

Let $\mu:S(G)\longrightarrow S(G)/H$ be the natural map. Since
$S(G)/H= S(\widehat{G})/\widehat{H}$ we have $\mu
(e_1)=\widehat{\mu}(e_1)=\widehat{\mu}(e_2)=\mu(e_2)$, so there
exists $h\in H$ such that $e_1 = he_2$. Since $e_1$ and $e_2$ have
a common vertex $v$, we have $h\in H_v$. Therefore, $e_2 =
h^{-1}e_1 = h^{-1}\gamma e_2$ and so $\gamma_1:= h^{-1}\gamma\in
\overline{H_{e_2}}$.

Recall that stabilizers of vertices in $H$ are finitely generated
and so by Theorem 1 in \cite{B-C-K-98},  $I$ is finitely
generated. Since $G_v$ is subgroup separable (Proposition
\ref{distinguido sob conjugacao}), the groups $I$ and $H_{e_2}$
are closed in the profinite topology of $G_v$. Therefore by the
main result in \cite{C-2001} $IH_{e_2}$ is closed in the profinite
topology of $G_v$, i.e., $\overline{IH_{e_2}}\cap G_v = IH_{e_2}$.
Thus $h= \gamma\gamma_1^{-1}\in IH_{e_2}$ and so there exist
$h_1\in I$, $h_2\in H_{e_2}$ such that $h=h_1h_2$. Let
$\gamma^{+}=\gamma h_1^{-1}$, then $\gamma^{+}e_1= \gamma h_1^{-1}
e_1=\gamma h_2 h^{-1}e_1=\gamma e_2 = e_1$, and therefore
$\gamma^{+}\in \overline{H}_{e_1}$. We also have $\gamma^{+}=
\gamma h_1^{-1}\in \overline{I}$, and since $I$ and $H_{e_1}$ are
finitely generated in $H_v$, by Proposition 2.4 in \cite{W-Z} for
the cases $d=1,2,7,11$ and by Lemma \ref{intersection} for the
case $d=3$ we have $\gamma^{+}\in \overline{\underset{e\in
P^{+}}{\cap}H_e}$. Therefore replacing $\gamma$ with $\gamma^{+}$
and $g_2$ with $h_1g_2h_1^{-1}$ we may assume that $\gamma \in
\overline{\underset{e\in P^{+}}{\cap}H_e}$.

\medskip

Thus we always can assume that $P$ has sufficiently large length.
Let $P$ be a geodesic of $T_{g_1}$ that contains $e$ and $g_1e$.
Let $I = \underset{e\in P}{\bigcap}H_e$ and consider $D = I\cap
g_1Ig_1^{-1}$. Note that $g_1Ig_1^{-1}$ is the intersection of the
stabilizers of the path $g_1P$. Now  by what was done above, we
can suppose that $\gamma\in \overline{D}$ (we use $P^{+} = P\cup
g_1P$ and if $e'=g_1e$, then $H_{e'}= g_1H_eg_1^{-1}$). Observe
that $I$ is cyclic, because the intersection of two edge groups in
$S(G)$ is cyclic, by Lemma 4.1 in \cite{W-Z}. Now we prove that
$g_1$ normalizes $D$. Indeed, since $D$ is closed in $G$, then $D
= \underset{N \unlhd_f G}{\cap}DN$. Thus it is enough to proof
that $DN$ is normalized by $g_1$, for each $N \unlhd_f G$. Indeed,
let $N\unlhd_f G$ and consider the quotient $G/N$. Then the group
$DN/N$ have the same index $m$ in  the groups $IN/N$ and
$I^{g_1^{-1}}N/N$, because these groups are conjugate. If $xN$
generates $IN/N$, then $(xN)^m$ and $(g_1N)(xN)^m(g_1N)^{-1}$
generate $DN/N$ and we conclude that $DN$ is normalized by $g_1$.

Let $d$ be a generator of $D$ and write $E = \langle d,
g_1\rangle$. We proof that $E$ is abelian. Indeed, if $d$ does not
centralize $g_1$,  $g_1dg_1^{-1}= d^{-1}$ and so $E$ is a
generalized dihedral group. However, by Corollary \ref{diedral}
the Bianchi groups do not contain subgroups isomorphic to the
generalized dihedral group, a contradiction.

Therefore $d$ centralizes  $g_1$ and consequently $E$ is abelian.
This implies that $\overline{E}$ is abelian and so $g_1 = g_2$, as
needed.
\end{proof}

 \begin{thm}\label{principal}
 Let $G$ be a torsion free group commensurable with a Euclidean  Bianchi
 group. Suppose that $G$ does not contain subgroups isomorphic to the
 generalized dihedral group. Then  $G$ is conjugacy separable.
 \end{thm}
 \begin{proof} Since $G$ is commensurable with $\Gamma_d$ (note that commensurability
 with $\Gamma_3$
 is equivalent to commensurability with $\Gamma_{3,2}$),  there exist isomorphic subgroups
  $M$, $N$ of finte index in $G$ and $\Gamma_d$ respectively. Then
  by the preceding proposition $N$ is conjugacy separable. Consequently $M$ is conjugacy
  separable. Then $M\leq G$ satisfies all hypothesis of Theorem \ref{G s.c.}, hence
 $G$ is conjugacy separable.
  \end{proof}
  \begin{thm}\label{principal2}
  Let $G\leq PGL_2(\cc)$ be a torsion free group commensurable with $\Gamma_d$,
  $d=1,2, 3,7,11$. Then
  $G$ is conjugacy separable.
  \end{thm}
\begin{proof} Follows from Theorem \ref{principal} and Corollary
\ref{diedral}.
 \end{proof}
\begin{cor}
The figure eight knot group is conjugacy separable.
\end{cor}
\begin{proof}
The figure eight knot group is a torsion free subgroup of index 12
in $\Gamma_3=PSL_2(O_3)$ (see page 60 in \cite{M-R}), so by
Theorem \ref{principal2} $K$ is conjugacy separable.
\end{proof}

\begin{rem} In Theorem \ref{principal2} $PGL_2(\cc)$ can be replaced by $GL_2(\cc)$ because
a subgroup of $GL_2(\cc)$ containing the generalized dihedral
group $\langle x,y\mid x^2=y^2\rangle$ must have $x^2\in
Z(GL_2(\cc)$ by Lemma \ref{lema0} and so can not be commensurable
with $\Gamma_d$. Indeed, every subgroup of finite index in
$\Gamma_d$ has finite center.
\end{rem}

\bigskip
\noindent

\begin{minipage}[u]{0.40\linewidth}
  Department of Mathematics\\
 University of Amazonas\\
 Avenida General Rodrigo Oct\'avio Jord\~ao Ramos, 3000\\
 Aleixo 69000-077\\
 Manaus-AM\\
 Brazil\\
 sheilachagas@ufam.edu.br
    \end{minipage}\hfill
    \begin{minipage}[u]{0.40\linewidth}
    \noindent Department of Mathematics\\
 University of Bras\'{\i}lia\\
Brasilia-DF 70910-900\\
 Brazil\\
pz@mat.unb.br
    \end{minipage}


\begin{thebibliography}{99}
\bibitem[B-71]{B} Burns, R. G., {\it On Finitely Generated Subgroups of Free Products},
 J. Austral. Math. Soc., {\bf 12} (1971) 358-364.
\bibitem[B-C-K-98]{B-C-K-98} Burns R. G., Chau E. C., Kam S.-M., {\it On The Rank of Intersection
of Subgroups of a Free Product of Groups}, J. Pure Appl. Algebra,
{\bf 124} (1998),  31-45.
\bibitem[C-2001]{C-2001} Coulbois T., {\it Free Product, Profinite Topology and
 Partial Action of Groups on Relational
Structures: A Connection Between Model Theory and Profinite
Topology}, World Scientific, Singapore 2002, 349-361.

\bibitem[D-80]{Dicks} Dicks W., {\it Groups, Trees and Projective
Modules}, Lecture Notes in Mathematics,  No. 790, Springer-Verlag,
1980.
\bibitem[D-D-89]{Dicks2} Dicks W., Dunwoody M. J., {\it Groups
Acting on Graphs}, Cambrige Univ. Press., Cambrige, 1989.
\bibitem[Dy-79]{Dyer-79} Dyer J. L., {\it Separating Conjugates
in Amalgamated Free by Finite Groups}, J. London Math. Soc., 2,
{\bf 20} (1979) 215-222.
\bibitem[F-89]{Fine} Fine B., {\it Algebraic Theory of the Bianchi
Groups}, Marcel Dekker, 1989.
\bibitem[G-M-S-2003]{G-M-S} Gitik R., Margolis S.W., Steinberg
B., {\it On the Kurosh theorem and separability properties},
Journal of Pure and Applied Algebra, {\bf 179} (2003) 87-97.
\bibitem[H-R-85]{H-R} Herfort W. H., Ribes L., {\it Torsion Elements and Centralizers in Free
Products of Profinite Groups}, J. Reine Angew. Math., {\bf 358}
(1985) 155-161.
\bibitem[L-R-2004]{L-R} Long D. D., Reid A. W., {\it Subgroup Separability
and Virtual Retractions of Groups}, Preprint, 23pp.
\bibitem[M-58]{M} Mal'cev A. I., {\it On Homormorphisms onto Finite
Groups}, Uchen. Zap. Ivanovskogo Gos. Ped. Inst., {\bf 18} (1958)
40-60.
\bibitem[M-R-02]{M-R} Maclachlan C., Reid A. W., {\it The Arithmetic of
Hyperbolic $3$-Manifolds}, Graduate Texts in Mathematics, {\bf
219}, Springer, 2002.
\bibitem[N-92]{Niblo} Niblo G. A., {\it Separability Properties
of Free Groups and Surface Groups}, J. Pure Appl. Algebra, {\bf
78} (1992) 77-84.
\bibitem[R-71]{R-71} Remeslennikov V. N., {\it Groups that are Residually  Finite with
Respect to Conjugacy}, Siberian Math. J., {\bf 12} (1971) 783-792.
\bibitem[Ro-69]{Ro} Romanovskii V. N., {\it On the Residual Finiteness of Free Products
with Respect to Subgroups}, Izv. Akad. Nauk SSSR Ser. Mat., {\bf
33} (1969) 1324-1329.
\bibitem[R-T-V-95]{R-T-V-95} Raptis E., Talelli O., Varsos D.,
{\it On the Conjugacy Separability of Certain Graphs of Groups},
J. Algebra, {\bf 199} (1998) 327-336.
\bibitem[R-S-Z-98]{R-S-Z} Ribes L., Segal D., Zalesskii P. A., {\it Conjugacy
Separability And Free Products of Groups with Cyclic
Amalgamation}, J. London Math. Soc., 2, {\bf  57} (1998) 609-628.

\bibitem[RZ]{RZ} L. Ribes, P. Zalesskii,   {\it Profinite groups.}
Ergebnisse der Mathematik und ihrer Grenzgebiete, 3. Folge {\bf 40},
Springer-Verlag, Berlin, (2000).

\bibitem[R-Z-96]{R-Z} Ribes L., Zalesskii P. A., {\it Conjugacy
Separability of Almagamated Free Products of Groups}, Journal of
Algebra, {\bf 179} (1996) 751-774.
\bibitem[R-Z-93]{R-Z-93} Ribes L., Zalesskii P. A., {\it On the
Profinite Topology of a Free Group}, Bull. London Math. Soc., {\bf
25} (1993) 37-43.
\bibitem[R-Z-04]{R-Z-04}  Ribes L., Zalesskii P.A., {\it Profinite topologies in
Free Products of Groups}. International Conference on Semigroups
and Groups in honor of the 65th birthday of Prof. John Rhodes.
Internat. J. Algebra Comput.  14  (2004),  no. 5-6, 751--772.
\bibitem[S-70]{Serre-70} Serre J.P., {\it Le Probleme de Groupes de Congruence pour
$SL_2$}, Ann. of Math., {\bf 92} (1970) 489-657.
\bibitem[S-77]{Serre} Serre J.P., {\it Trees}, Springer
Monographs in Mathematics, Springer, 1977.
\bibitem[Sc-03]{Scarth} Scarth M. R., {\it PhD Thesis: Normal
Congruence Subgroups of Bianchi Groups and Related Groups},
Faculty of Science, University of Glasgow (2003).
\bibitem[W-Z-98]{W-Z} Wilson J. S., Zalesskii P. A., {\it Conjugacy Separability of
Certain Bianchi Groups and HNN Extensinos}, Math. Proc. Camb.
Phil. Soc., {\bf 123} (1998) 227-242.
\bibitem[Z-M-89]{Z-M-89} Zalesskii P. A., Mel'nikov O. V., {\it Subgroups
of Profinite Groups Acting on Trees}, Math. Sb., {\bf 63} (1989)
405-424.
\bibitem[Z-M-90]{Z-M-90} Zalesskii P. A., Mel'nikov O. V., {\it Fundamental Groups of Graphs
of Profinite Groups}, Algebra i Analiz, {\bf 1} (1989) 117-134.
\bibitem[Z-T-95]{Z-T}  Zalesskii P. A., Tavgen O. I., {\it Closed
Orbits and Finite Approximability With Respect to Conjugacy of
Free Almagamated Products}, Mathematical Notes, {\bf 58} (1995)
1042-1048.
\bibitem[Z-92]{Z-92} Zalesskii P. A., {\it Open Subgroups of Free
Profinite Products}, Contemporary Mathematics, {\bf 131} (1992)
473-491.
\bibitem[Z-89]{Z-89} Zalesskii P. A., {\it A Geometrical Characterization of Free
 Constructions of Profinite Groups}, Siber. Math. J., {\bf 30} (1989)
 73-84.
\end{thebibliography}
\end{document}